\documentclass[12pt, a4paper]{article}

\usepackage{amsfonts}
\usepackage{amsmath}
\usepackage{amssymb}
\usepackage[english]{babel}
\usepackage{amsthm}

\newtheorem {theorem}{Theorem}[section]
\newtheorem {lemma}[theorem]{Lemma}
\newtheorem {proposition}[theorem]{Proposition}
\newtheorem {corollary}[theorem]{Corollary}

\theoremstyle{definition}
\newtheorem* {remark}{Remark}

\title{Asymptotics for some combinatorial characteristics
of the convex hull of a Poisson point process in the Clifford torus}

\author{Alexander Magazinov \thanks{Supported by the Russian government project 11.G34.31.0053 and RFBR grant 11-01-00633-a.}}

\begin{document}

\maketitle

\begin{abstract}
Let $\mathcal P_\lambda$ be a homogeneous Poisson point process of rate $\lambda$ in the Clifford torus $T^2\subset \mathbb E^d$.
Let $(f_0, f_1, f_2, f_3)$ be the $f$-vector of $conv\,\mathcal P_\lambda$ and let $\bar{v}$
be the mean valence of a vertex of the convex hull. Asymptotic expressions for $\mathsf E\, f_1$, $\mathsf E\, f_2$, $\mathsf E\, f_3$ and
$\mathsf E\, \bar{v}$ as $\lambda\to\infty$ are proved in this paper.

\end{abstract}

\pagestyle{myheadings}
\thispagestyle{empty}

\section{Introduction}

Recently Poisson-Voronoi tessellations became an object for extensive investigations. The first non-trivial result concerning 
Poisson-Voronoi tessellations belongs to J.~L.~Meijering. His paper \cite{mei} shows that a typical
cell of a Poisson-Voronoi tessellation of the 3-dimensional Euclidean space $\mathbb E^3$ has an expectation of number of facets 
equal to 
$$\frac{48\pi^2}{35}+2 = 15.5354\ldots\;. $$

A survey \cite{wwi} contains a number of further results concerning Poisson-Voronoi tessellations.

It is possible to consider Voronoi tessellations of a sphere or a hyperbolic space of constant curvature as well as 
Voronoi tessellations of a Euclidean space. Given a locally finite set $A$ in a sphere, Euclidean space or a hyperbolic
space of constant curvature, it is possible to consider an associated Delaunay triangulation. It is known (see, for example, \cite{rog})
that the following statements are equivalent.

\begin{enumerate}
\item A subset $B\subset A$ spans a face of Delaunay triangulation. 
\item A set of points equidistant to all points of $B$ contains a face of Voronoi tessellation associated with $A$. 
\end {enumerate}

Therefore the notions of Voronoi tessellation and Delaunay triangulation are dual to each other.

Consider a finite set $A$ of points in general position in the sphere
$$S^d_r=\{(\xi_1, \xi_2, \ldots, \xi_{d+1})\subset \mathbb E^{d+1}:\xi_1^2+\xi_2^2+\ldots+\xi_{d+1}^2=r^2\}.$$
A subset $B\subset A$ determines a face of the Delaunay triangulation associated with $A$ if and only if 
$conv\, B$ is a face of $conv\, A$.

N.~Dolbilin and M.~Tanemura (see \cite{dta}) studied convex hulls of finite subsets of the Clifford torus $T^2$ embedded in
$\mathbb E^4$. Since $T^2\subset S^3_{\sqrt{2}}$, this case can be considered as an additional restriction for a finite
subset of $S^3_{\sqrt{2}}$ generating the Delaunay triangulation (or the Voronoi tessellation). For a special class of point sets in $T^2$ called {\it regular sets} \cite{dta} completely describes the combinatorial structure 
of the convex hull. 

In addition, the convex hull of the Poisson point process within $T^2$ has been explored by numeric methods.
Dolbilin and Tanemura considered the average number $\bar{f}$ of 2-faces of a cell in the corresponding Voronoi tessellation of
$\mathbb S^3$, which is exactly the average degree of a vertex in the Delaunay triangulation. A strong linear relation between
$\bar{f}$ and $\log_{10} N$ ($N = 4\pi^2\lambda$ is the average number of points) has been observed, and the obtained regression
formula was
$$\bar{f} \approx -2.419308 + 9.971915 \log_{10} N.$$

In other words, the simulation has shown that the mean valence of a vertex of the convex hull (or the mean number of hyperfaces
of a Poisson-Voronoi cell) is likely to have an expectation $O^*(\ln \lambda)$, as the rate of the process $\lambda$ tends to infinity.

Here and further $F_1=O^*(F_2)$ means
that $\limsup\limits_{\lambda\to\infty} 
\max\left ( \left |\frac {F_1}{F_2} \right|, \left |\frac {F_2}{F_1} \right| \right ) <\infty$.

N.~Dolbilin suggested the author to prove the conjecture on the logarithmic growth of the mean valence of a vertex.

In this paper this conjecture and several related results are proved.

\section{Notation and main results}

In the 4-dimensional Euclidean space $\mathbb E^4$ consider the two-dimensional Clifford torus
$$T^2=\{(\cos \phi, \sin \phi, \cos \psi, \sin \psi):-\pi<\phi, \psi\leq \pi\}.$$ 
Clearly, $T^2$ is a submanifold of the three-dimensional sphere 
$$S^3_{\sqrt{2}}=\{(\xi_1,\xi_2,\xi_3,\xi_4):\xi_1^2+\xi_2^2+\xi_3^2+\xi_4^2=2\}.$$

$T^2$ has a locally Euclidean planar metric and, consequently, the natural Borel measure $mes_2$, where $mes_2(T^2)=4\pi^2$. 

Consider a random point set $\mathcal P \subset T^2$. For every Borel-measurable set $A\subset T^2$ define a random
variable
$$n(A) = n_{\mathcal P}(A) = |\mathcal P \cap A|.$$ 

Denote by $Pois(\nu)$ the Poisson distribution with rate parameter $\nu$, i.e. the distribution of a random variable
$\zeta_\nu$ such that 
$$\mathsf P(\zeta_\nu = j) = e^{-\nu} \frac{\nu^j}{j!}\quad \text{for} \quad j=0,1,2,\ldots .$$
 
Say that $\mathcal P = \mathcal P_{\lambda}$ is the (homogeneous) Poisson point process of rate $\lambda>0$ if the random 
variable $n(A)$ is distributed according to $Pois\bigl ( \lambda\, mes_2(A)\bigr ) $ law for every Borel-measurable set $A\subset T^2$. 

Call a polytope in $\mathbb E^4$ {\it generic} if it is a simplicial 4-polytope, or a simplex of dimension at most 3,
or an empty polytope. Remind the notion of $f$-vector of a 4-polytope and extend it to the cases of other generic
polytopes.

The $f$-vector of a 4-polytope $P$ is a 4-vector $(f_0, f_1, f_2, f_3)$, where $f_i$ is the number of $i$-faces of $P$
for $i=0,1,2,3$. By definition, for the 3-dimensional simplex put the $f$-vector equal to $(4,6,4,2)$, for the 2-dimensional 
simplex --- $(3,3,1,0)$, for the segment --- $(2,1,0,0)$, for one point --- $(1,0,0,0)$, and for empty polytope --- $(0,0,0,0)$.

If $P = conv\, \mathcal P_{\lambda}$, then $P$ is almost surely a generic polytope, and therefore $(f_0, f_1, f_2, f_3)$ is
a well-defined random vector.

Call the event $n(T^2)\leq 4$ a {\it degenerate case} and the complementary event $n(T^2)>4$, respectively,
a {\it non-degenerate case}.

\begin{remark} The reason to choose $f_3=2$ for a 3-dimensional simplex is the convenience to treat it
as a polytope with 2 hyperfaces equal to this simplex. The other components were chosen to satisfy Dehn-Sommerville
equations (see, for example, \cite{bpa}). The $f$-vectors for other polytopes occuring in degenerate cases were chosen rather arbitrarily according to the common idea of simplices in dimensions lower than 3.
\end{remark}

The main results of this paper are below.

\begin{theorem}
\label{t1}
The number of hyperfaces of 
$conv\, \mathcal P_{\lambda}$ has a magnitude of expectation $O^*(\lambda\ln\lambda)$ as $\lambda$ tends to infinity.
\end{theorem}

\begin{theorem}
\label{t2}
The numbers of 1-faces and 2-faces 
of $conv\, \mathcal P_{\lambda}$ both have magnitudes of expectation $O^*(\lambda\ln\lambda)$ as $\lambda$ tends to infinity.
\end{theorem}

In addition, one can easily observe that the value of $f_0$ (i.e. the number of vertices) for the polytope 
$conv\, \mathcal P_{\lambda}$ is exactly $n(T^2)$. Therefore
$$\mathsf E\, f_0 = \mathsf E\, n(T^2) = 4\lambda\pi^2,$$
as $n(T^2)$ is $Pois(4\lambda\pi^2)$-distributed (expectations of Poisson random variables are computed, for example, in \cite{shi}).

\begin{remark} For a random polytope $conv\, \mathcal P_{\lambda}$ the asymptotics of the expectation of $f$-vector as
$\lambda \to \infty$ is now completely described.
\end{remark}

The other combinatorial characteristiñ of a polytope is mean valence of its vertices. More precisely, given
a polytope $P$ in $\mathbb E^4$ (possibly, empty) with $f$-vector $(f_0, f_1, f_2, f_3)$, consider the value 
$$\bar{v}= \bar{v}(P) = \left\{
\begin{array}{ll}
\frac{2f_1}{f_0},\quad & \text{if}\quad f_0\neq 0,\\
0,\quad & \text{if}\quad f_0=0.
\end{array}
\right.
$$  
Then $\bar{v}$ is called {\it the mean valence of vertex} of $P$. If $P = conv\, \mathcal P_{\lambda}$, then
$\bar{v}=\bar{v}(conv\, \mathcal P_{\lambda})$ is a random variable.

\begin{theorem}
\label{t3}
The expectation of the mean 
valence of a vertex of $conv\, \mathcal P_{\lambda}$ has asymptotics $\mathsf E\,\bar{v}=O^*(\ln\lambda)$ as $\lambda$ tends to infinity.
\end{theorem}

\begin{remark} Theorem \ref{t3} provides an answer to the problem proposed by Dolbilin and Tanemura.
\end{remark}

Here and further the designations of all combinatorial characteristics apply to the random polytope $conv\, \mathcal P_{\lambda}$.

\section{Integral expressions for $\mathsf E\,f_3$ and $\mathsf E\,\bar{v}$}

Let $(T^2)^4$ be the fourth Cartesian power of $T^2$ with natural measure $mes_8$. Let $X\subset (T^2)^4$ be the set of all points 
$x=(x_1, x_2, x_3, x_4)$, 
where $x_i\in T^2$ such that points $x_1, x_2, x_3, x_4$ are affinely independent in $\mathbb E^4$. 

For every $x\in X$
denote by $p(x)$ a hyperplane spanned by points $x_1$, $x_2$, $x_3$, $x_4$. 
It is obvious that $X$ is open in $(T^2)^4$. Moreover, it is easily seen that $(T^2)^4\setminus X$ has a zero
measure. 

Denote by $\Pi^+(x)$ and $\Pi^-(x)$ the two half-spaces determined by $p(x)$ for every $x\in X$.

The sets 
$$C^+(x) = T^2\cap \Pi^+(x)\quad \text{and}\quad C^-(x) = T^2\cap \Pi^-(x)$$
are called {\it caps}. 

Without loss of generality, assume that for every $x\in X$
$$mes_2(C^+(x))\leq mes_2(C^-(x)).$$

Let $G:X\to\mathbb R$ be a function determined by 
$$G(x) = mes_2(C^+(x)).$$ 
Clearly, $G(x)$ is continuous on $X$.  

The integral expressions for $E\,f_3$ and $E\,\bar{v}$ will be obtained by using the famous Slivnyak-Mecke formula. This
formula was proved for the first time in \cite{mec}, and in \cite{dst} it is stated as follows.

\begin{proposition}[{\rm Slivnyak-Mecke formula}] 
\label{sm}
Let $\mathcal X$ be a space with measure $\mu$. Suppose $\mathcal N_{\mathcal X}$ is a space of all locally finite 
point configurations in $\mathcal X$. Consider a Poisson point process $\mathcal P_\mu$ within $\mathcal X$ corresponding 
to the measure $\mu$. Then for every measurable function 
$F: \mathcal X^s \times \mathcal N_{\mathcal X} \to [0, \infty)$ holds
\begin{equation}\label{01}
\begin{split}
\mathsf E\, \sum\limits_{\{x_1, x_2, \ldots, x_s\}\subset \mathcal P_\mu}^{\neq} 
F\left(x_1, x_2, \ldots, x_s, \mathcal P_\mu\setminus \{x_1, x_2, \ldots, x_s\} \right)=\\
=\int\limits_{\mathcal X^s} \mathsf E\, \bigl( F\left(x_1, x_2, \ldots, x_s, \mathcal P_\mu \right)\bigr) \, d\mu(x_1)d\mu(x_2)\ldots d\mu(x_s).
\end{split}
\end{equation}
\end{proposition}

The sign $\neq$ here stands for summation over all $s$-tuples of {\it distinct} points.

Throughout the proofs of Lemma \ref{l1} and Lemma \ref{l2} $\lambda$ is assumed to be a fixed positive real number.

\begin{lemma}
\label{l1}
\begin{equation} \label{14}
\mathsf E\,f_3= \frac {1}{24} \int\limits_{(T^2)^4} \lambda^4\left(e^{-\lambda G(x)}+e^{-\lambda(4\pi^2-G(x))}\right)\, dx.
\end{equation} 
\end{lemma}

\begin{proof}
Apply the Slivnyak-Mecke formula (\ref{01}) for 
$$\mathcal X = T^2,\quad s=4, \quad \mu=\lambda\cdot mes_2 \quad \text{and}$$
$$F(x_1, x_2, x_3, x_4, X) = \mathbf{1}_{C^+(x)\cap X\subset \partial C^+(x)}+\mathbf{1}_{C^-(x)\cap X\subset \partial C^-(x)}.$$

If $n(T^2)>4$ and $x_1, x_2, x_3, x_4$ are distinct points of $\mathcal P_\lambda$ then it is not hard to see that almost surely
\begin{multline*}
F\left(x_1, x_2, x_3, x_4, \mathcal P_\nu\setminus \{x_1, x_2, x_3, x_4\} \right)=\\
=\left \{
\begin{array}{ll}
1,\quad & \text{if $x_1, x_2, x_3, x_4$ span a hyperface of $\mathcal P_\lambda$} \\
0,\quad & \text{otherwise}.
\end{array} \right.
\end{multline*}

If $n(T^2)=4$ and $x_1, x_2, x_3, x_4$ are the four points of $\mathcal P_\lambda$ then almost surely
$$F\left(x_1, x_2, x_3, x_4, \mathcal P_\nu\setminus \{x_1, x_2, x_3, x_4\} \right)=2.$$

Finally, if $n(T^2)<4$ then there are no quadruples in $\mathcal P_\lambda$ and the left part of (\ref{01}) is
an empty sum.

Therefore in every case 
\begin{equation}\label{02}
\sum\limits_{\{x_1, x_2, x_3, x_4\}\subset \mathcal P_\lambda}^{\neq} 
F\left(x_1, x_2, x_3, x_4, \mathcal P_\lambda\setminus \{x_1, x_2, x_3, x_4\} \right)= 24f_3, 
\end{equation}
since the quadruple $(x_1, x_2, x_3, x_4)$ can be ordered in 24 different ways.

Moreover, by definition of a Poisson point process,
\begin{equation}\label{03}
\begin{split}
\mathsf E\, \mathbf{1}_{C^+(x)\cap X\subset \partial C^+(x)} = e^{-\lambda\cdot mes_2(C^+(x))}& = e^{-\lambda G(x)},\\
\mathsf E\, \mathbf{1}_{C^-(x)\cap X\subset \partial C^-(x)} = e^{-\lambda\cdot mes_2(C^-(x))}& = e^{-\lambda(4\pi^2-G(x))}.
\end{split}
\end{equation}

Substitution of (\ref{02}) and (\ref{03}) into (\ref{01}) gives the statement of Lemma \ref{l1}.

\end{proof}

For $\nu>0$ let $\zeta_\nu$ be distributed as $Pois(\nu)$. Denote
\begin{equation}\label{4hs}
h(\nu)= \mathsf E\, \frac{1}{\zeta_\nu+4} = \sum\limits_{j=0}^{\infty} \frac{\nu^j}{j!(j+4)}e^{-\nu}.
\end{equation}

Direct computation of the sum in (\ref{4hs}) gives
\begin{equation}\label{4h}
h(\nu)=\frac{1}{\nu}-\frac{3}{\nu^2}+\frac{6}{\nu^3}-\frac{6-6e^{-\nu}}{\nu^4}.
\end{equation}

Obviously, $h(\nu)$ is continuous for $\nu>0$.

\begin{lemma}
\label{l2}
\begin{multline} \label{43}
\mathsf E\,\bar{v}
= \frac {1}{12} \int\limits_{(T^2)^4} \lambda^4\left(e^{-\lambda G(x)}h\bigl(4\lambda\pi^2-\lambda G(x)\bigr)+
e^{-4\lambda\pi^2+\lambda G(x)}h\bigl(\lambda G(x)\bigr)\right)\, dx \; +\\
2 - \mathsf P\bigl(n(T^2)=2\bigr) - 2\mathsf P\bigl(n(T^2)<2\bigr).
\end{multline} 
\end{lemma}

\begin{proof}

The Dehn-Sommerville equations \cite[Section 1.2]{bpa} hold for $conv \, \mathcal P_\lambda$ almost 
surely in non-degenerate cases as well as in the case $n(T^2)=4$. These equations imply $f_1=f_3+f_0$. Therefore in these cases
\begin{equation}\label{111}
\bar{v}=2\frac{f_3}{f_0}+2.
\end{equation}

Apply the Slivnyak-Mecke formula (\ref{01}) for 
$$\mathcal X = T^2,\quad s=4, \quad \mu=\lambda\cdot mes_2 \quad \text{and}$$
$$F(x_1, x_2, x_3, x_4, X) = \left(\mathbf{1}_{C^+(x)\cap X\subset \partial C^+(x)}+\mathbf{1}_{C^-(x)\cap X\subset \partial C^-(x)}
\right)\cdot \frac {1}{|X\cup \{x_1, x_2, x_3, x_4\}|}.$$
The substitution gives the following identity
\begin{multline}\label{46}
24 \mathsf E\left (\frac{f_3}{f_0} \mid n(T^2)\geq 4\right) \cdot P\bigl(n(T^2)\geq 4\bigr) = \\
\int\limits_{(T^2)^4} \lambda^4\left(e^{-\lambda G(x)}h\bigl(4\lambda\pi^2-\lambda G(x)\bigr)+
e^{-4\lambda\pi^2+\lambda G(x)}h\bigl(\lambda G(x)\bigr)\right)\, dx.
\end{multline}

According to the law of total probability,
\begin{multline}\label{47}
\mathsf E \, \bar{v} = 2\mathsf E\left (\frac{f_3}{f_0} \mid n(T^2)\geq 4\right)\cdot \mathsf P\bigl(n(T^2)\geq 4\bigr)+ 2 \mathsf P\bigl(n(T^2)\geq 4\bigr)+\\
2 \mathsf P\bigl(n(T^2)=3 \bigr) + \mathsf P\bigl(n(T^2)=2\bigr).
\end{multline}

By (\ref{46}), the first summand at the right side of (\ref{47}) is equal to the integral in (\ref {43}). Further,
$$\mathsf P\bigl(n(T^2)<2 \bigr) + \mathsf P\bigl(n(T^2)=2 \bigr) + \mathsf P\bigl(n(T^2)=3 \bigr)+ \mathsf P\bigl(n(T^2)\geq 4 \bigr) =1,$$
therefore
\begin{multline*}
2 \mathsf P\bigl(n(T^2)\geq 4\bigr)+  2 \mathsf P\bigl(n(T^2)=3 \bigr) + \mathsf P\bigl(n(T^2)=2\bigr) =\\
2 - \mathsf P\bigl(n(T^2)=2\bigr) - 2 \mathsf P\bigl(n(T^2)<2\bigr),
\end{multline*}
so the remaining parts of the right sides of  (\ref {43}) and (\ref {47}) are equal as well. 

\end{proof}

\section{Estimates for the measure function}

To proceed we need two statements about the caps. Both of them are proved by a fairly simple computation, so the proofs
are given in the Appendix.

\begin{lemma}\label{l3} 

The following statements hold:

\begin{enumerate}

\item[\rm 1.] For every cap $C^+(x)$ (respectively,  $C^-(x)$)
there exist $a, b\geq 0$, $\phi_0, \psi_0$  satisfying $a^2+b^2\geq 2$ and
$-\pi<\phi_0,\psi_0\leq \pi$ such that 
$$C^+(x)=\left\{
(\phi, \psi)\in T^2: a^2 \sin^2 \frac{\phi - \phi_0}{2} + b^2 \sin^2 \frac{\phi - \phi_0}{2} \leq 1\right\}, $$
and, respectively,
$$C^-(x)=\left\{
(\phi, \psi)\in T^2: a^2 \sin^2 \frac{\phi - \phi_0}{2} + b^2 \sin^2 \frac{\phi - \phi_0}{2} \geq 1\right\}.$$ 
where $\gamma_1,\gamma_2>0$ and do not depend on $x$.

\item[\rm 2.] For every $a, b\geq 0$, $\phi_0, \psi_0$  satisfying $a^2+b^2\geq 2$ and
$-\pi<\phi_0,\psi_0\leq \pi$ the sets
$$\left\{(\phi, \psi)\in T^2: a^2 \sin^2 \frac{\phi - \phi_0}{2} + b^2 \sin^2 \frac{\phi - \phi_0}{2} \leq 1\right\} \quad \text{and} $$
$$\left\{(\phi, \psi)\in T^2: a^2 \sin^2 \frac{\phi - \phi_0}{2} + b^2 \sin^2 \frac{\phi - \phi_0}{2} \geq 1\right\}$$ 
are caps.

\end{enumerate}

\end{lemma}

\begin{remark}
$a$, $b$, $\phi_0$, $\psi_0$ can be now considered as functions $a(x)$, $b(x)$, $\phi_0(x)$, 
$\psi_0(x)$ of the argument $x\in X$.
\end{remark}

\begin{lemma}\label{l4} 

There exist positive constants $\gamma_1, \gamma_2$ such that for every $x\in X$ holds
$$\gamma_1<(a(x)+1)(b(x)+1)G(x)<\gamma_2.$$

\end{lemma}

For every $t\in \mathbb R$ define
$$M(t)=mes_8\{x\in X: G(x)<t\},$$
$$N(t)=mes_8\{x\in X: G(x)<t\; \text{and}\; \min(a(x), b(x))<100\},$$
$$L(t)=mes_8\{x\in X: G(x)<t\; \text{and}\; \min(a(x), b(x))\geq 100\}.$$

It is easily seen that $M(t)=N(t)=L(t)=0$ for $t<0$ and $M(t)=N(t)+L(t)$ for every $t\in \mathbb R$.

The main goal of the present section is to estimate $M(t)$. We estimate $N(t)$ and $L(t)$ separately in Lemma \ref{l5} and Lemma \ref{l6}.

\begin{lemma}
\label{l5}
There exists $\gamma_3>0$ such that
$$N(t)<\gamma_3 t^3$$
for every $0<t<\frac{1}{2}$.
\end{lemma}

\begin{lemma}
\label{l6}
There exist $\gamma_4, \gamma_5>0$ such that
$$\gamma_4 t^3|\ln t|< L(t)< \gamma_5 t^3|\ln t|$$
for every $0<t<\frac{1}{2}$.
\end{lemma}

Before the proofs we give an estimate of $M(t)$ as a corollary.

\begin{corollary}

There exist positive constants $\gamma_6, \gamma_7$ such that
$$\gamma_6 t^3|\ln t|<M(t)< \gamma_7 t^3|\ln t|$$
for every $0<t<\frac{1}{2}$.

\end{corollary}

\begin{proof}[Proof of Lemma \ref{l5}]
Introduce the functions
$$N_1(t)=mes_8\{x\in X: G(x)<t\; \text{and}\; a(x)<100\},$$
$$N_2(t)=mes_8\{x\in X: G(x)<t\; \text{and}\; b(x)<100\}.$$

Obviously, $N_1(t)=N_2(t)$ and $N(t)\leq N_1(t)+N_2(t)$.

Suppose 
$$0<t\leq \frac{\gamma_1}{1000\pi}.$$

Let $a(x)<100$ and $G(x)<t$. Lemma \ref{l4} implies  
$$b(x)\geq \frac {\gamma_1}{(a+1)G(x)}-1 > \frac {\gamma_1}{200t}.$$

By Lemma \ref{l3}, cap $C^+(x)$ is described by the inequality
$$a(x)^2 \sin^2 \frac{\phi-\phi_0}{2}+ b(x)^2 \sin^2 \frac{\psi-\psi_0}{2}\leq 1.$$

From the last inequality follows that every point of $C^+(x)$ with coordinates $(\phi, \psi)$ satisfies
$$\left| \sin \frac{\psi-\psi_0}{2} \right| \leq \frac {1}{b} < \frac {200t}{\gamma_1}.$$

Hence
$$C^+(x)\subset \left\{ (\phi, \psi)\in T^2:
   \left| \sin \frac{\psi-\psi_0}{2} \right|< \frac {200t}{\gamma_1} \right\} = S(t,x).$$

A set $S\subset T^2$ is called a {\it strip} if there exist $\psi_1\in (-\pi, \pi]$ and $d\in (-1, 1)$
such that 
$$S=\{(\phi, \psi)\in T^2: \cos (\psi-\psi_1) \leq d.\}$$
The {\it centerline} of $S$ is the line
$\psi=\psi_1$,
and $2\arccos d$ is the {\it width} of $S$.

$S(t,x)$ is obviously a strip of width
$$w(t)=4\arcsin \frac {200t}{\gamma_1}< \frac {1000\pi t}{\gamma_1}$$
and the centerline of $S(t,x)$ is described by the equation $\psi=\psi_0$.

Let 
$$k=k(t)=\lceil \frac{2\pi}{w(t)} \rceil.$$

Consider $k$ strips $S_1, S_2, \ldots S_k \subset T^2$ of width
$2w(t)$ each such that $S_j$ has centerline $\psi=-\pi+\frac {2\pi j}{k}$. 

It is obvious that $S(t,x)\subset S_j$, where $j$ is the nearest integer to 
$\frac{k (\psi_0+\pi)}{2\pi}$ and $S_0=S_k$.

Let $x=(x_1, x_2, x_3, x_4)$, where $x_i\in T^2$. Obviously, every $x_i\in \partial C^+(x)$, therefore $x\in S_j^4$.

Finally, 
\begin{multline*}
N_1(t)=mes_8\{x\in X: G(x)<t\; \text{and}\; a(x)<100\}\leq \\
mes_8 \left(\bigcup\limits_{j=1}^k(t) S_j^4 \right) \leq k(t)(4\pi w(t))^4 \leq (4\pi)^5w(t)^3\leq \gamma'_3 t^3.
\end{multline*} 

Then
$$N(t)\leq 2N_1(t) \leq 2\gamma'_3 t^3.$$

The case 
$$0<t\leq \frac{\gamma_1}{1000\pi}$$
is proved completely.

Suppose 
$$\frac{\gamma_1}{1000\pi}<t<\frac{1}{2}.$$

Obviously, 
$$N(t)\leq mes_8\bigl((T^2)^4 \bigr)= 256 \pi^8.$$ 

Then
$$N(t)<  256 \pi^8 \left( \frac{1000\pi}{\gamma_1} \right)^3 t^3, $$
and Lemma \ref{l5} is now proved completely.

\end{proof}

\begin{proof}[Proof of Lemma \ref{l6}]
Suppose $\min(a(x), b(x))\geq 100$. Assume 
$$x=(x_1, x_2, x_3, x_4),\quad \text{where}\quad x_i=(\phi_i, \psi_i)\in T^2 \quad
\text{for} \quad i=1,2,3,4.$$

Let 
$$\alpha(x)=\frac {1}{a(x)}, \quad \beta(x)=\frac{1}{b(x)}.$$

By assumptions, $0<\alpha(x),\beta(x)<\frac{1}{100}$.

Lemma \ref{l4} easily implies that there exist $\gamma'_1, \gamma'_2>0$ such that
\begin{equation}\label{31}
\gamma'_1\alpha(x)\beta(x)<G(x)<\gamma'_2\alpha(x)\beta(x).
\end{equation}

Since $x_i = (\phi_i, \psi_i) \in \partial C^+(x)$ for $i=1,2,3,4$, then
$$\frac{\sin^2 \frac{\phi_i-\phi_0}{2}}{\alpha^2}+\frac{\sin^2 \frac{\psi_i-\psi_0}{2}}{\beta^2}=1.$$

Therefore we can define parameters $-\pi < \theta_i \leq \pi$ for $i=1,2,3,4$ such that
$$\sin \frac {\phi_i-\phi_0}{2}=\alpha\cos\theta_i\quad \text{and} \quad
\sin \frac {\psi_i-\psi_0}{2}=\beta\sin\theta_i.$$

It is not hard to see that every point $x\in X$ parametrized by 8 numbers 
$$(\alpha, \beta, \phi_0, \psi_0, \theta_1, \theta_2, \theta_3, \theta_4)$$
can be uniquely parametrized by another 8 numbers
$$(\phi_1, \psi_1, \phi_2, \psi_2, \phi_3, \psi_3, \phi_4, \psi_4).$$

Since there are two parametrizations of a point $x\in X$, consider a Jacobi matrix between
these parametrizations. The elements are computed as follows:
$$\frac{\partial \phi_i}{\partial \phi_0}=1,\;  \frac{\partial \psi_i}{\partial \phi_0}=0;$$
$$\frac{\partial \phi_i}{\partial \psi_0}=0,\;  \frac{\partial \psi_i}{\partial \psi_0}=1;$$
$$\frac{\partial \phi_i}{\partial \alpha}=\frac{2\cos\theta_i}{\cos\frac {\phi_i-\phi_0}{2}};\;
     \frac{\partial \psi_i}{\partial \alpha}=0;$$
$$\frac{\partial \phi_i}{\partial \beta}=0;\;
     \frac{\partial \psi_i}{\partial \beta}=\frac{2\sin\theta_i}{\cos\frac {\psi_i-\psi_0}{2}};$$
$$\frac{\partial \phi_i}{\partial \theta_i}=\frac{-2\alpha\sin\theta_i}{\cos\frac {\phi_i-\phi_0}{2}};\;
     \frac{\partial \psi_i}{\partial \theta_i}=\frac{2\beta\cos\theta_i}{\cos\frac {\psi_i-\psi_0}{2}};\;
          \frac{\partial \phi_i}{\partial \theta_j}=\frac{\partial \psi_i}{\partial \theta_j}=0.$$
          
Therefore
\begin{multline*}
J=\left|\frac{D(\phi_1, \psi_1, \phi_2, \psi_2, \phi_3, \psi_3, \phi_4, \psi_4)}
{D(\phi_0, \psi_0, \alpha, \beta, \theta_1, \theta_2, \theta_3, \theta_4)}\right|=\\
\left|
\begin{smallmatrix}
1&0&1&0&1&0&1&0\\
0&1&0&1&0&1&0&1\\
\frac{2\cos\theta_1}{\cos\frac {\phi_1-\phi_0}{2}}&0&
   \frac{2\cos\theta_2}{\cos\frac {\phi_2-\phi_0}{2}}&0&
      \frac{2\cos\theta_3}{\cos\frac {\phi_3-\phi_0}{2}}&0&
         \frac{2\cos\theta_4}{\cos\frac {\phi_4-\phi_0}{2}}&0\\
0&\frac{2\sin\theta_1}{\cos\frac {\psi_1-\psi_0}{2}}&
   0&\frac{2\sin\theta_2}{\cos\frac {\psi_2-\psi_0}{2}}&
      0&\frac{2\sin\theta_3}{\cos\frac {\psi_3-\psi_0}{2}}&
         0&\frac{2\sin\theta_4}{\cos\frac {\psi_4-\psi_0}{2}}\\
\frac{-2\alpha\sin\theta_1}{\cos\frac {\phi_1-\phi_0}{2}}& \frac{2\beta\cos\theta_1}{\cos\frac {\psi_1-\psi_0}{2}}&
   0&0&0&0&0&0\\
0&0&
   \frac{-2\alpha\sin\theta_2}{\cos\frac {\phi_2-\phi_0}{2}}& \frac{2\beta\cos\theta_2}{\cos\frac {\psi_2-\psi_0}{2}}&
      0&0&0&0\\
0&0&0&0&
   \frac{-2\alpha\sin\theta_3}{\cos\frac {\phi_3-\phi_0}{2}}& \frac{2\beta\cos\theta_3}{\cos\frac {\psi_3-\psi_0}{2}}&
      0&0\\
0&0&0&0&0&0&
   \frac{-2\alpha\sin\theta_4}{\cos\frac {\phi_4-\phi_0}{2}}& \frac{2\beta\cos\theta_4}{\cos\frac {\psi_4-\psi_0}{2}}\\
\end{smallmatrix}
\right|.
\end{multline*}

Direct computation shows that
\begin{multline*}
J=\sum\limits_{(i\,j\,k\,l)} 64\, \mathrm{sign}\,(i\,j\,k\,l)\, \alpha^2\beta^2 \cdot 
\frac {1}{
                  \prod\limits_{m=1}^4 \cos\frac{\phi_m-\phi_0}{2} \cos\frac{\psi_m-\psi_0}{2} } \times\\
\cos^2 \theta_i \cos \theta_j \sin^2 \theta_k \sin \theta_l \cos\frac{\phi_j-\phi_0}{2}\cos\frac{\psi_l-\psi_0}{2},
\end{multline*}
where $(i\,j\,k\,l)$ runs through all permutations of $(1\, 2\, 3\, 4)$.

From (\ref{31}) easily follows that
\begin{multline*}
mes_8\left\{x\in (T^2)^4: \alpha(x)\beta(x)<\frac{t}{\gamma'_2} \; \text{and} \max(\alpha(x), \beta(x))<\frac{1}{100} \right\}\leq 
L(t) \leq \\
mes_8\left\{x\in (T^2)^4: \alpha(x)\beta(x)<\frac{t}{\gamma'_1} \; \text{and} \max(\alpha(x), \beta(x))<\frac{1}{100} \right\}.
\end{multline*}

Therefore
\begin{multline*}
\int\limits_{\substack
   {\max(\alpha, \beta)<\frac{1}{100}\\
    \alpha\beta<\frac{t}{\gamma'_2} }} d\phi_1 d\psi_1 d\phi_2 d\psi_2 d\phi_3 d\psi_3 d\phi_4 d\psi_4 \leq L(t) \leq\\
\int\limits_{\substack
   {\max(\alpha, \beta)<\frac{1}{100}\\
    \alpha\beta<\frac{t}{\gamma'_1} }} d\phi_1 d\psi_1 d\phi_2 d\psi_2 d\phi_3 d\psi_3 d\phi_4 d\psi_4.
\end{multline*}  

In variables 
$(\alpha, \beta, \phi_0, \psi_0, \theta_1, \theta_2, \theta_3, \theta_4)$
the last inequality can be written as follows
\begin{multline*}
\int\limits_{\substack
   {\max(\alpha, \beta)<\frac{1}{100}\\
    \alpha\beta<\frac{t}{\gamma'_2} }} 
   \int\limits_{\substack
   {\phi_0, \psi_0 \in (-\pi, \pi]\\
    \theta_{1,2,3,4} \in (-\pi, \pi]}} |J|\, d\alpha d\beta d\phi_0 d\psi_0 d\theta_1 d\theta_2 d\theta_3 d\theta_4 \leq L(t) \leq\\
\int\limits_{\substack
   {\max(\alpha, \beta)<\frac{1}{100}\\
    \alpha\beta<\frac{t}{\gamma'_1} }} 
   \int\limits_{\substack
   {\phi_0, \psi_0 \in (-\pi, \pi]\\
    \theta_{1,2,3,4} \in (-\pi, \pi]}} |J|\, d\alpha d\beta d\phi_0 d\psi_0 d\theta_1 d\theta_2 d\theta_3 d\theta_4.
\end{multline*}

Let 
\begin{multline} \label{32}
J_1=\frac{J}{\alpha^2\beta^2}=\sum\limits_{(i\,j\,k\,l)} 64\, \mathrm{sign}\,(i\,j\,k\,l)\, \cdot 
\frac {1}{
                  \prod\limits_{m=1}^4 \cos\frac{\phi_m-\phi_0}{2} \cos\frac{\psi_m-\psi_0}{2} } \times\\
\cos^2 \theta_i \cos \theta_j \sin^2 \theta_k \sin \theta_l \cos\frac{\phi_j-\phi_0}{2}\cos\frac{\psi_l-\psi_0}{2}.
\end{multline}
Then $J_1$ can be considered as a function 
$J_1(\alpha, \beta, \phi_0, \psi_0, \theta_1, \theta_2, \theta_3, \theta_4)$.

Obviously, with $\gamma = \gamma'_1$ or $\gamma'_2$
\begin{multline}\label{321}
\int\limits_{\substack
   {\max(\alpha, \beta)<\frac{1}{100}\\
    \alpha\beta<\frac{t}{\gamma} }} 
   \int\limits_{\substack
   {\phi_0, \psi_0 \in (-\pi, \pi]\\
    \theta_{1,2,3,4} \in (-\pi, \pi]}} |J|\, d\alpha d\beta d\phi_0 d\psi_0 d\theta_1 d\theta_2 d\theta_3 d\theta_4 = \\
\int\limits_{\substack
   {\max(\alpha, \beta)<\frac{1}{100}\\
    \alpha\beta<\frac{t}{\gamma} }} \alpha^2\beta^2\, d\alpha d\beta 
   \int\limits_{\substack
   {\phi_0, \psi_0 \in (-\pi, \pi]\\
    \theta_{1,2,3,4} \in (-\pi, \pi]}} |J_1|\, d\phi_0 d\psi_0 d\theta_1 d\theta_2 d\theta_3 d\theta_4.
\end{multline}

Since $\max(\alpha, \beta)<\frac{1}{100}$, then 
$$|\sin \frac{\phi_m-\phi_0}{2}|<\frac {1}{100}\quad \text{and} \quad
|\sin \frac{\psi_m-\psi_0}{2}|<\frac {1}{100}$$ 
for $m=1,2,3,4$. 

Hence  
$$\cos \frac{\phi_m-\phi_0}{2}> \frac{4999}{5000} \quad \text{and} \quad
\cos \frac{\psi_m-\psi_0}{2}> \frac{4999}{5000}.$$

Consequently,
\begin{multline*}
\left |\cos^2 \theta_i \cos \theta_j \sin^2 \theta_k \sin \theta_l \cos\frac{\phi_j-\phi_0}{2}\cos\frac{\psi_l-\psi_0}{2} - \right.\\
\left. \phantom{\frac{\phi_j-\phi_0}{2}} \cos^2 \theta_i \cos \theta_j \sin^2 \theta_k \sin \theta_l \right| \leq \frac{2}{5000}.
\end{multline*}

Applying (\ref{32}), we obtain the following sequence of inequalities which are independent from $\alpha, \beta$:
\begin{multline*}
J_1\geq \\
64 \left| \sum\limits_{(i\,j\,k\,l)} \, \mathrm{sign}\,(i\,j\,k\,l)\, \cdot 
\cos^2 \theta_i \cos \theta_j \sin^2 \theta_k \sin \theta_l \cos\frac{\phi_j-\phi_0}{2}\cos\frac{\psi_l-\psi_0}{2} \right| \geq\\
64 \left| \sum\limits_{(i\,j\,k\,l)} \, \mathrm{sign}\,(i\,j\,k\,l)\, \cdot 
\cos^2 \theta_i \cos \theta_j \sin^2 \theta_k \sin \theta_l \right| - 64\cdot 24\cdot\frac{2}{5000}.
\end{multline*}

Let $\theta_m=(m-2)\pi/2$. Then
$$\left| \sum\limits_{(i\,j\,k\,l)} \, \mathrm{sign}\,(i\,j\,k\,l)\, \cdot 
\cos^2 \theta_i \cos \theta_j \sin^2 \theta_k \sin \theta_l \right| = 4.$$

Consequently, there exists some neighbourhood of point $(-\frac{\pi}{2}, 0, \frac{\pi}{2}, \pi)$ in
coordinates $(\theta_1, \theta_2, \theta_3, \theta_4)$ such that 
$|J_1|>64$ in this neighbourhood, and the neighbourhood is independent from $\alpha, \beta, \phi_0, \psi_0$. 

Therefore there exist positive constants $\gamma'_4, \gamma'_5$ independent of $\alpha, \beta$ and satisfying 
\begin{equation}\label{323}
\gamma'_4 <
\int\limits_{\substack
   {\phi_0, \psi_0 \in (-\pi, \pi]\\
    \theta_{1,2,3,4} \in (-\pi, \pi]}} |J_1|\, d\phi_0 d\psi_0 d\theta_1 d\theta_2 d\theta_3 d\theta_4 < \gamma'_5
\end{equation}
for every $0<\alpha, \beta < \frac{1}{100}$.

Inequalities (\ref{321}) and (\ref{323}) together imply
\begin{equation}\label{33}
 \gamma'_4\cdot \int\limits_{\substack
   {\max(\alpha, \beta)<\frac{1}{100}\\
    \alpha\beta<\frac{t}{\gamma'_2} }} \alpha^2\beta^2\, d\alpha d\beta < L(t) <
    \gamma'_5\cdot \int\limits_{\substack
   {\max(\alpha, \beta)<\frac{1}{100}\\
    \alpha\beta<\frac{t}{\gamma'_1} }} \alpha^2\beta^2\, d\alpha d\beta. 
\end{equation}
    
If $\tau<\frac{1}{10000}$ then
$$\int\limits_{\substack
   {\max(\alpha, \beta)<\frac{1}{100}\\
    \alpha\beta<\tau }} \alpha^2\beta^2\, d\alpha d\beta = \frac{\tau^3}{9}-\frac{2 \ln 100}{9} \tau^3 + \frac{1}{9}\tau^3|\ln \tau|.$$
    
Consequently, as $t\to 0$, the main terms in left and right parts of (\ref{33}) have order $t^3|\ln t|$. Therefore
$\frac{L(t)}{t^3|\ln t|}$ is bounded from above and below in some interval $(0, \varepsilon)$ by two positive constants.

In the segment $[\varepsilon, \frac 12]$ the functions $L(t)$ and $t^3|\ln t|$ are continuous and positive. Consequently,
the quotient $\frac{L(t)}{t^3|\ln t|}$ in this segment is also bounded from above and below by two positive constants (but,
probably, not the same as in the previous paragraph). However, combining the cases $t \in (0, \varepsilon)$ and
$t \in [\varepsilon, \frac 12]$ allows to conclude that $\frac{L(t)}{t^3|\ln t|}$ is bounded from above and below
by some positive constants in the whole interval $(0, \frac 12]$. Thus Lemma \ref{l6} is proved.

\end{proof}

\section{Proofs of main results} Now proceed with the proofs of Theorems \ref{t1}, \ref{t2} and \ref{t3}.

\begin{proof}[Proof of Theorem \ref{t1}] 

It is obvious that
\begin{multline*}
\int\limits_{(T^2)^4} \lambda^4 e^{-\lambda G(x)}\, dx<
\int\limits_{(T^2)^4} \lambda^4\left(e^{-\lambda G(x)}+e^{-\lambda(4\pi^2-\lambda G(x))}\right)\, dx<\\
2\int\limits_{(T^2)^4} \lambda^4 e^{-\lambda G(x)}\, dx.$$
\end{multline*}

Then, by Lemma \ref{l1},
$$\mathsf E\, f_3 = O^*\left(  \int\limits_{(T^2)^4} \lambda^4 e^{-\lambda G(x)}\, dx  \right).$$

According to \cite{gra}, the identity
$$\int\limits_{(T^2)^4} \lambda^4 e^{-\lambda G(x)}\, dx=
\int\limits_{\mathbb R} \lambda^4 e^{-\lambda t}\, dM(t),$$
holds, where the right side is a Stieltjes integral.

Since $G(x)$ is the measure of $C^+(x)$, then $0<G(x)\leq 2\pi^2$ holds for every $x\in X$. 
Therefore $M(t)$ is a constant for $t\leq 0$ and for $t\geq 2\pi^2$. Hence
$$\int\limits_{(T^2)^4} \lambda^4 e^{-\lambda G(x)}\, dx=\int\limits_{0}^{2\pi^2} \lambda^4 e^{-\lambda t}\, dM(t).$$

Thus
\begin{equation}\label{56}
\mathsf E\, f_3 = O^* \left( \int\limits_{0}^{2\pi^2} \lambda^4 e^{-\lambda t}\, dM(t) \right).
\end{equation}

Since $M(t)$ is non-decreasing, $e^{-\lambda t}$ is decreasing and continuous, then integration by parts is possible
and gives
\begin{equation}\label{57}
\int\limits_{0}^{2\pi^2} \lambda^4 e^{-\lambda t}\, dM(t) =
\lambda^4 M(2\pi^2)e^{-2\lambda \pi^2} + \lambda^5 \int\limits_{\frac{1}{2}}^{2\pi^2} e^{-\lambda t}M(t)\, dt
+ \lambda^5 \int\limits_{0}^{\frac{1}{2}} e^{-\lambda t}M(t)\, dt.
\end{equation}

Obviously, as $\lambda \to \infty$,
\begin{equation}\label{58}
\lambda^4 M(2\pi^2)e^{-2\lambda \pi^2}=o(1),\quad
\lambda^5 \int\limits_{\frac{1}{2}}^{2\pi^2} e^{-\lambda t}M(t)\, dt =o(1),
\end{equation}

$$\lambda^5 \int\limits_{0}^{\frac{1}{2}} e^{-\lambda t}M(t)\, dt = 
O^*\left( \lambda^5\int\limits_{0}^{\frac{1}{2}} e^{-\lambda t}t^3|\ln t|\, dt \right).$$

Let $u=e^{-\lambda t}$, then 
$$t=- \frac{\ln u}{\lambda}\quad \text{and} \quad dt=-\frac{du}{\lambda u}.$$

Therefore
$$\int\limits_{0}^{\frac{1}{2}} e^{-\lambda t}t^3|\ln t|\, dt = 
\int\limits_{e^{\frac{-\lambda}{2}}}^{1} \frac{|\ln^3 u|}{\lambda^4} \cdot (\ln\lambda -\ln(-\ln u))\, du=
O^*(\lambda^{-4}\ln\lambda).$$

Hence
\begin{equation}\label{59}
\lambda^5 \int\limits_{0}^{\frac{1}{2}} e^{-\lambda t}M(t)\, dt = O^*(\lambda \ln\lambda).
\end{equation}

Substitution of (\ref{58}) and (\ref{59}) into (\ref{57}) gives
$$\int\limits_{0}^{2\pi^2} \lambda^4 e^{-\lambda t}\, dM(t) = O^*(\lambda \ln\lambda).$$

Thus, according to (\ref{56}),
$$\mathsf E\, f_3 = O^*(\lambda \ln\lambda),$$
which is the statement of Theorem \ref{t1}.

\end{proof}

\begin{proof}[Proof of Theorem \ref{t2}] 

From Dehn-Sommerville equations for a simplicial 4-polytope follows that
$$f_2=2f_3+r_2 \quad \text{and} \quad
f_1=f_3-f_0+r_1,$$
where random variables $r_1$ and $r_2$ are errors of degenerate cases, i.e. $r_1=r_2=0$ almost surely if $n(T^2)>4$ 
and  $r_1, r_2<10$ almost surely. Since 
$$\lim\limits_{\lambda\to\infty}\mathsf P(n(T^2)\leq 4)=0,$$ 
then 
$$\mathsf E\, r_1 = o(1), \quad \mathsf E\, r_2 = o(1),\quad
    \mathsf E\, f_3 = O^*(\lambda \ln\lambda)$$
as $\lambda\to\infty$ by Theorem \ref{t1}. Also,
$$\mathsf E\, f_0 = \mathsf E\, n(T^2) = 4\lambda\pi^2$$
since $n(T^2)$ is distributed as $Pois(4\lambda\pi^2)$.

Finally,
$$\mathsf E\, f_2 = 2\mathsf E\, f_3 + \mathsf E\, r_2 = O^*(\lambda \ln\lambda),$$
$$\mathsf E\, f_1= \mathsf E\, f_3 - \mathsf E\, f_0 + \mathsf E\, r_1 = O^*(\lambda \ln\lambda),$$
and Theorem \ref{t2} is proved.

\end{proof}

\begin{proof}[Proof of Theorem \ref{t3}]
Notice that
$$h\bigl(\lambda G(x)\bigr)<\frac{1}{4} \quad \text{and} \quad
e^{-4\lambda\pi^2+\lambda G(x)}\leq e^{-2\lambda\pi^2}.$$

These inequalities imply the estimate
$$ \int\limits_{(T^2)^4} e^{-4\lambda\pi^2+\lambda G(x)}h\bigl(\lambda G(x)\bigr)\, dx \leq 64\pi^8 e^{-2\lambda\pi^2} = o(1)$$
as $\lambda\to\infty$.

Further,
$$2\lambda\pi^2\leq 4\lambda\pi^2-\lambda G(x) < 4\lambda\pi^2.$$

Therefore from (\ref{4h}) follows 
$$h\bigl(4\lambda\pi^2-\lambda G(x)\bigr)=O^*(\lambda^{-1}).$$

Consequently,
$$ \int\limits_{(T^2)^4} \lambda^4 e^{-\lambda G(x)}h\bigl(4\lambda\pi^2-\lambda G(x)\bigr)  \, dx =
O^*\left( \lambda^3 \int\limits_{(T^2)^4} e^{-\lambda G(x)} \, dx \right) = O^*(\ln \lambda),$$
as the integral in the middle part was estimated in the proof of Theorem \ref{t1}.

Obviously, 
$$\mathsf P\bigl(n(T^2)=2\bigr)=o(1), \quad \mathsf P\bigl(n(T^2)<2\bigr)=o(1).$$

Now Theorem \ref{t3} easily follows from (\ref{43}) because every summand in this identity was estimated.

\end{proof}

\section*{Acknowledgements} Author is grateful to N.~Dolbilin who stated the problem and gave many useful
recommendations on the text, to V.~Dolnikov for a constant and enormous help during the work on this paper,
to I.~B\'ar\'any and M.~Tanemura for useful discussions.

\section*{Appendix. Structure of caps} 

This section is devoted to obtaining analytical description and measure estimates for the caps.

\begin{proof}[Proof of Lemma \ref{l3}] 

Suppose
$$p(x)=\left\{
(\xi_1, \xi_2, \xi_3, \xi_4)\in\mathbb E^4 : a_1 \xi_1+a_2 \xi_2+b_1 \xi_3+b_2 \xi_4=c \right\}$$
where $c\geq 0$.

Then
$$\partial C^+(x)=\partial C^-(x)= \left\{ (\phi, \psi) \in T^2:
a_1\cos \phi+ a_2 \sin\phi + b_1 \cos\psi + b_2\sin\psi = c \right\}.$$

The equation for $\partial C^+(x)$ can be rewritten as
$$a'\cos(\phi-\phi_0)+b'\cos(\psi-\psi_0)=c,$$
where $a'=\sqrt{a_1^2+a_2^2}$ and $b'=\sqrt{b_1^2+b_2^2}$.

Since 
$$\cos(\phi-\phi_0)=1-2\sin^2 \frac {\phi-\phi_0}{2} \quad \text{and} 
\cos(\psi-\psi_0)=1-2\sin^2 \frac {\psi-\psi_0}{2},$$ 
the previous equation is eqivalent to
$$a'\sin^2\frac {\phi-\phi_0}{2}+b'\sin^2\frac {\psi-\psi_0}{2}=\frac {a'+b'-c}{2}.$$

The set $\partial C^+(x)$ contains infinitely many points, therefore
$$0\leq c< a'+b'.$$

If $c=0$ then $p(x)$ passes through
the origin and therefore divides $T^2$ into equal parts.  Consequently, 
\begin{multline*}
mes_2\left (\left\{
  (\phi, \psi)\in T^2: a'\sin^2\frac {\phi-\phi_0}{2}+b'\sin^2\frac {\psi-\psi_0}{2}\leq \frac {a'+b'}{2} \right\} \right) = \\
mes_2\left (\left\{
  (\phi, \psi)\in T^2: a'\sin^2\frac {\phi-\phi_0}{2}+b'\sin^2\frac {\psi-\psi_0}{2}\geq \frac {a'+b'}{2}\right\} \right).
\end{multline*}

Therefore for $c>0$
\begin{multline*}
mes_2\left (\left\{
  (\phi, \psi)\in T^2: a'\sin^2\frac {\phi-\phi_0}{2}+b'\sin^2\frac {\psi-\psi_0}{2}\leq \frac {a'+b'-c}{2} \right\} \right) < \\
mes_2\left (\left\{
  (\phi, \psi)\in T^2: a'\sin^2\frac {\phi-\phi_0}{2}+b'\sin^2\frac {\psi-\psi_0}{2}\geq \frac {a'+b'-c}{2}\right\} \right).
\end{multline*}

Since $mes_2(C^+(x))\leq mes_2(C^-(x))$,
$$C^+(x)=\left\{
(\phi, \psi)\in T^2: a'\sin^2\frac {\phi-\phi_0}{2}+b'\sin^2\frac {\psi-\psi_0}{2}\leq \frac {a'+b'-c}{2} \right\}.$$

Let 
$$a=\sqrt{\frac{2a'}{a'+b'-c}}\quad \text{and}\quad  b=\sqrt{\frac{2b'}{a'+b'-c}}.$$ 

Then
$$C^+(x)=\left\{(\phi, \psi)\in T^2 : a^2 \sin^2 \frac{\phi - \phi_0}{2} + b^2 \sin^2 \frac{\psi - \psi_0}{2} \leq 1\right\}$$
and, respectively,
$$C^-(x)=\left\{(\phi, \psi)\in T^2 : a^2 \sin^2 \frac{\phi - \phi_0}{2} + b^2 \sin^2 \frac{\psi - \psi_0}{2} \geq 1\right\},$$
hence statement 1 of Lemma \ref{l3}.

All the computations are obviously invertible, and performing them in the inverse order gives statement 2 of Lemma \ref{l3}.

\end{proof}

\begin{proof}[Proof of Lemma \ref{l4}] 

Without loss of generality assume $\phi_0=\psi_0=0$.

Consider the case $a=0$. Then
$$C^+(x)=\left\{(\phi, \psi)\in T^2: b^2 \sin^2 \frac{\phi}{2} \leq 1\right\},$$
or, equivalently,
$$C^+(x)=\left\{(\phi, \psi)\in T^2: |\psi|\leq 2\arcsin \frac{1}{b} \right\}.$$

Consequently,
$$(a+1)(b+1)G(x)=8\pi(b+1)\arcsin \frac {1}{b}.$$

Since $\frac {1}{b}<\arcsin\frac {1}{b}< \frac {\pi}{2} \cdot \frac {1}{b}$ and $b\geq \sqrt{2}$,
$$8\pi<(a+1)(b+1)G(x)<4\pi^2 \left(1+\frac {\sqrt{2}}{2} \right), $$
and the case $a=0$ is completely proved. The case $b=0$ is similar.

Now suppose $a>0$ and $b>0$. Since 
$$ \frac {|\phi|}{\pi} \leq \left| \sin \frac {\phi}{2} \right| \leq \frac {|\phi|}{2}
\quad \text{and} \quad
\frac {|\psi|}{\pi} \leq \left| \sin \frac {\psi}{2} \right| \leq \frac {|\psi|}{2},$$
then the following inclusions hold:
$$C^+(x) \subset 
\left\{(\phi, \psi)\in T^2: 
a^2 \left(\frac {\phi}{2} \right)^2\leq \frac{1}{2} \text{ and } b^2 \left(\frac {\psi}{2} \right)^2\leq \frac{1}{2}
\right\},$$
$$C^+(x) \supset 
\left\{(\phi, \psi)\in T^2: 
a^2 \left(\frac {\phi}{\pi} \right)^2\leq 1 \text{ and } b^2 \left(\frac {\psi}{\pi} \right)^2\leq 1
\right\}.$$

Therefore
$$ \min \left( 2\pi, \frac{2\sqrt{2}}{a} \right) \cdot \min \left( 2\pi, \frac{2\sqrt{2}}{b} \right)  
\leq G(x) \leq
\min \left( 2\pi, \frac{2\pi}{a} \right) \cdot \min \left( 2\pi, \frac{2\pi}{b} \right).$$

It is easy to check that
$$ \min \left( 2\pi, \frac{2\sqrt{2}}{a} \right)= \frac {1}{\max \left( \frac{1}{2\pi}, \frac{a}{2\sqrt{2}} \right)}\geq 
\frac {1}{ \frac{1}{2\pi}+\frac{a}{2\sqrt{2}} }= \frac{2\pi\sqrt{2}}{\pi a + \sqrt{2}},$$
$$ \min \left( 2\pi, \frac{2\pi}{a} \right) \leq
\frac {2}{ \frac{1}{2\pi}+\frac{a}{2\pi} }= \frac{4\pi}{a+1},$$
and, similarly,
$$ \min \left( 2\pi, \frac{2\sqrt{2}}{b} \right)\geq \frac{2\pi\sqrt{2}}{\pi a + \sqrt{2}},$$
$$\min \left( 2\pi, \frac{\pi\sqrt{2}}{b} \right) \leq \frac{4\pi}{b+1}.$$

Finally,
$$
8\leq \frac{2\pi\sqrt{2}(a+1)}{\pi a + \sqrt{2}} \cdot \frac{2\pi\sqrt{2}(b+1)}{\pi b + \sqrt{2}}\leq (a+1)(b+1)G(x) 
\leq 16\pi^2,
$$
which completes the proof of Lemma \ref{l4}.

\end{proof}

\end{document}